\documentclass{amsart}
\usepackage{amsmath}
\usepackage{amsfonts}
\usepackage{amssymb}
\usepackage{graphicx}
\usepackage[mathscr]{eucal}
\newtheorem{theorem}{Theorem}
\theoremstyle{plain}

\newtheorem{lemma}{Lemma}

\numberwithin{equation}{section}
\begin{document}
\title[Log-concavity of $f$-vectors of ordinary polytopes]{Log-concavity of face vectors of cyclic and ordinary polytopes}
\author{L\'{a}szl\'{o} Major}\address[]{L\'{a}szl\'{o} Major \newline\indent Institute of Mathematics \newline\indent Tampere University of Technology \newline\indent PL 553, 33101 Tampere, Finland}\email[]{laszlo.major@tut.fi}
\date{Aug 6, 2011}

\begin{abstract}
Ordinary polytopes are known as a non-simplicial generalization of the cyclic polytopes. The face vectors of ordinary polytopes are shown to be log-concave.
\end{abstract}
\maketitle
In the late 1950's, Motzkin (according to Bj\"orner \cite{bjo}) conjectured that  the face vectors of convex polytopes are unimodal. This conjecture was apparently disproved by Danzer (presented in a lecture in Graz in 1964, see \cite{zie}), consequently the question raised: ''Which natural classes of polytopes still have unimodal $f$-vectors?'' For an overview of unimodality and log-concavity questions related to convex polytopes see e.g. \cite{bjo, eck, wer, zie}.

The special shape of the $f$-vectors of cyclic polytopes makes them
applicable to certain constructions to present non-unimodal convex polytopes \cite{zieg}.
Nevertheless, the $f$-vectors of cyclic polytopes themselves are unimodal and log-concave. It seems natural to ask, whether the ordinary polytopes, the generalizations of the cyclic polytopes introduced by Bisztriczky \cite{bis}, have log-concave $f$-vectors. 

For any $d$-polytope $P$, we denote its $f$-vector by $f(P)$ whose $j^{th}$ component $f_j(P)$ is the number
of faces of dimension $j$ in $P$ for $j=-1,0,\ldots,d-1$. We have $f_{-1}(P)=1$ for the improper face $\emptyset$. The vector $(f_{-1},f_0,\ldots,f_{d-1}) $ is called \emph{log-concave} if $f_{i-1}f_{i+1}\leq f_i^2$ for all $-1<i<d-1$. In dimension less than 5, the $f$-vector of any polytope is log-concave (\cite{wer}), therefore from now on we assume that $d\geq 5$.
According to the concept of the $h$-vector of simplicial polytopes, we define the \textit{$h$-vector} for any polytope by the same linear relations:
 $$h_i:=\sum_{j=0}^i(-1)^{i-j}\binom{d-j}{d-i}f_{j-1}, \hspace{2mm}\text{for } \hspace{2mm}0\leq i \leq  d.$$
 It is easy to invert the above equations to return to 
{the components of the $f$-vector:} $$f_{j}=\sum_{i=0}^d \binom{d-i}{d-j-1}h_i, \hspace{2mm}\text{for } \hspace{2mm}-1\leq j \leq  d-1.$$
In other words, to know the $h$-vector is equivalent to know the $f$-vector. These relations can be easily made visible by using the following method based on Stanley's technique (see \cite{sta}).
We replace the first $d+1$ bordering 1's on the right-hand side of Pascal's triangle by the components of the $h$-vector, while the internal entries (indicated by $\triangle$ in the following figure) obey the usual "sum of two entries above" rule. The $f$-vector emerges in the $(d+1)^{th}$ row of this modified Pascal's tri\vspace{4mm}angle.

\setlength{\unitlength}{0.65mm}
\begin{picture}(100,61)(-8,1.5)

\put(90,58){$\textbf{h}$} 

\put(85,50){$1$}\put(95,50){$\textbf{v}$} 

\put(80,42){$1$}\put(90,42){$\triangle$} \put(100,42){$\textbf{e}$}

\put(75,34){$1$} \put(85,34){$\triangle$}\put(94,34){$\triangle$} \put(104.3,34){$\textbf{c}$}

\put(70,26){$1$} \put(80,26){$\triangle$} \put(90,26){$\triangle$} \put(100,26){$\triangle$} \put(110,26){$\textbf{t}$} 

\put(65,18){$1$} \put(75,18){$\triangle$} \put(85,18){$\triangle$}
\put(95,18){$\triangle$} \put(105,18){$\triangle$}\put(115,18){$\textbf{o}$}

\put(60,10){$1$}\put(70,10){$\triangle$} \put(79,10){$\triangle$}
\put(89.5,10){$\triangle$} \put(100,10){$\triangle$}\put(110,10){$\triangle$} \put(120,10){$\textbf{r}$}

 \put(55,2){$\textbf{f}$} \put(66,2){$\textbf{v}$}\put(75,2){$\textbf{e}$} \put(85,2){$\textbf{c}$}\put(96,2){$\textbf{t}$} \put(106,2){$\textbf{o}$}
\put(116,2){$\textbf{r}$}

\put(59,7.9){\line(1,0){64.2}}
\put(88.0,58.3){\line(3,-5){33.3}}
\end{picture}
\vspace{4mm}\\
This version of Stanley's technique is used by Lee in \cite{lee}. To give this triangle array  more formally we need some operators on vectors. Let us define the operator $N$ from $\mathbb{R}^s\times\mathbb{R}=\mathbb{R}^{s+1}$ to $\mathbb{R}^{s+1}$ for any $s$ by
\begin{equation*}\begin{split}
 N(\textbf{a},b)&=N((a_{-1},a_0,a_1,\ldots,a_{s-2}),b)=(x_{-1},x_0,x_1,\ldots,x_{s-1})\\
&=(a_{-1},a_{-1}+a_0,a_0+a_1,\ldots,a_{s-3}+a_{s-2},a_{s-2}+b).\\
\end{split}\end{equation*}
The operator $N$ provides a row of the modified Pascal's triangle using the previous row.
Let $\textbf{b}=(b_0,b_1,\ldots,b_r)\in \mathbb{R}^{r+1}$, then we define $F(\textbf{b})$ as the $(r+1)^{th}$ (last) vector in the vector sequence $$\textbf{b}(1)=(b_0),\textbf{b}(2)=N(\textbf{b}(1),b_1),\ldots,\textbf{b}(i+1)=N(\textbf{b}(i),b_i)),\ldots,\textbf{b}(r+1).$$ 
In this notation for the $h$-vector $\textbf{h}=(h_0,h_1,\ldots,h_d)\in \mathbb{R}^{d+1}$ of a $d$-polytope, Stanley \cite{sta} says that $F(\textbf{h})$ is the $f$-vector of the polytope.

Let $\textbf{a}=(a_{-1},a_0,a_1,\ldots,a_{s-2})\in \mathbb{R}^{s}$ and $\textbf{b}=(b_0,b_1,\ldots,b_r)\in \mathbb{R}^{r+1}$, then define 
$T(\textbf{a},\textbf{b})$ as the $(r+1)^{th}$ (last) vector in the vector sequence $$\textbf{b}(1)=N(\textbf{a},b_0),\textbf{b}(2)=N(\textbf{b}(1),b_1),\ldots,\textbf{b}(i+1)=N(\textbf{b}(i),b_i)),\ldots,\textbf{b}(r+1).$$ 
Again in the light of \cite{sta} for any $0<i<d$  the $f$-vector of a $d$-polytope is 
\begin{equation}\label{eq3}
T\big(F(h_0,h_1,\ldots,h_i),(h_{i+1},\ldots,h_d)\big).
\end{equation}

 Now, we recall some combinatorial results related to ordinary and cyclic polytopes, for further geometrical description see e.g. \cite{bay,bis,din}. For any ordinary $d$-polytope $P$ with $n+1$ vertices there exists a uniquely determined integer $k$ (the characteristic of $P$) such that the combinatorial structure and consequently the $f$-vector of $P$ depend only on the parameters $d,k$ and $n$. We denote the ordinary $d$-polytope $P$ with $n+1$ vertices and with characteristic $k$ by $P^{d,k,n}$. If $k=n$, then 
 $P^{d,k,n}$ is a cyclic $d$-polytope with $k+1$ vertices.
 \begin{theorem}[Bisztriczky \cite{bis}]\label{t1}Let $P^{d,k,n}$ be an ordinary polytope. If $d$ is even, then $P^{d,k,n}$ is  a cyclic polytope.\end{theorem}
Because of the above theorem, it is reasonable to deal separately with the even and odd dimensional cases. In the even dimensional case, we may take advantage of the following combinatorial result.
 \begin{theorem}[Brenti \cite{bre}]\label{t2}Let $\textbf{b}=(b_0,b_1,\ldots,b_d)$ be a positive log-concave vector. Then the vector $F(\textbf{b})$ is log-concave.
\end{theorem}
For a polytope $P$, the above statement says that the log-concavity of the $h$-vector  of $P$ implies the log-concavity of the $f$-vector of $P.$
The $h$-vector of the cyclic polytope $P^{d,n,n}$ is given as follows (see e.g.  \cite{zie})   
 \begin{equation}\label{eq1}h_{i}=\binom{n-d+i}{i}, \hspace{2mm}\text{for } \hspace{2mm}0\leq i \leq  \bigg\lfloor\frac d2\bigg\rfloor
 \end{equation}
 and for $\big\lfloor\frac d2\big\rfloor < i \leq d$ we have $h_i=h_{d-i}$ from the Dehn-Sommerville equations. 
Using the explicit formula $\binom{n}{k}=n!/k!(n-k)!$, it is easy to show that the $h$-vector of a cyclic polytope is log-concave. This fact together with Theorems \ref{t1} and \ref{t2} proves the following theorem.
\begin{theorem}\label{t3}If $d$ is even, then the $f$-vectors of ordinary $d$-polytopes are log-concave.
\end{theorem}
Now, we assume that $d=2m+1$ for some integer $m$. The above reasoning does not work in the odd dimensional case because the $h$-vectors of ordinary polytopes are not generally log-concave. In odd dimension, we need a specialization of a general result on triangular arrays of numbers proved by Kurtz in \cite{kur}:
\begin{theorem}\label{t4}Let $\textbf{a}$ be a log-concave vector and let $\textbf{b}=(b_0,b_1,\ldots,b_r)$ such that $b_0\geq b_1 \geq\ldots\geq b_r$. Then the vector $T(\textbf{a},\textbf{b})$ is log-concave.
\end{theorem} 
However, the above theorem is an immediate consequence of \cite{kur}, it can be easily proved directly by induction. 

Dinh computed the $f$-vectors of the ordinary polytopes in \cite{din}:
\begin{theorem}\label{t5}Let $n\geq k\geq d=2m+1>5$. The number of $j$-faces of the ordinary polytope $P^{n,k,d}$ is 
$$f_j(P^{n,k,d})=\phi_j(d,k)+(n-k)c_j(d,k),$$
where $\phi_j(d,k)$ denotes the number of $j$-faces of the cyclic $d$-polytope with $k+1$ vertices and  $c_j(d,k)=\sum^{m-1}_{i=0}\binom{m+1}{j-i}\binom{k-m-2}{i}$
\end{theorem} 
Originally, Dinh gave the numbers $c_j(d,k)$ in a more complicated  form, the above expression is calculated by Bayer \cite{bay}. Let $\textbf{u}(d,k)=(u_{-1},u_0,\ldots,u_{m-1})$ be defined by $u_i=\binom{k-m-2}{i}$ and let $\textbf{0}=(0,0,\ldots,0)$ be the zero vector with $m+1$ components. Then Bayer's formula can be expressed as follows.
\begin{equation}\label{eq2}\textbf{c}(d,k)=(c_{-1},c_0,\ldots,c_{d-1})=T(\textbf{u},\textbf{0}).
 \end{equation}
 Consider the example of the ordinary
polytope $P^{5,7,9}$. In this case $\textbf{u}=(0,1,3)$ and we compute $T(\textbf{u},\textbf{0})$ like this:

\setlength{\unitlength}{0.65mm}
\begin{picture}(100,41)(-8,1.5)

\put(75,34){$\big (0$} \put(85,34){$1$}\put(94,34){$3\big)$} \put(104.3,34){${{0}}$}

\put(70,26){$0$} \put(80,26){$1$} \put(90,26){$4$} \put(100,26){$3$} \put(110,26){${0}$} 

\put(65,18){$0$} \put(75,18){$1$} \put(85,18){$5$}
\put(95,18){$7$} \put(105,18){$3$}\put(115,18){${0}$}

\put(60,10){$\big(0$}\put(70,10){$1$} \put(79,10){$6$}
\put(89.5,10){$12$} \put(100,10){$10$}\put(110,10){$3\big)$} 

\put(35,10){$\textbf{c}(d,k)=$}
\put(59,34){$\textbf{u}=$}

\put(60,7.9){\line(1,0){56.2}}
\put(74,32.5){\line(1,0){36.3}}
\end{picture}

We turn our attention to the other quantity of the formula of Theorem \ref{t5}. According to (\ref{eq3}), for the $f$-vector of the cyclic polytope with $k+1$ vertices, we have

$$\textbf{f}(P^{d,k,k})=T\big(F(h_0,h_1,\ldots,h_m),(h_{m+1},\ldots,h_d)\big).$$
Let us define the vector $\textbf{v}(d,k)=(v_{-1},v_0,\ldots,v_{m-1})$ by $v_i=\binom{k-m}{i+1}$. Applying Stanley's technique to the first half of the $h$-vector given in (\ref{eq1}), we obtain that $F(h_0,h_1,\ldots,h_m)=\textbf{v}$. Therefore we have $\textbf{f}(P^{d,k,k})=T\big(\textbf{v},(h_{m+1},\vspace{-2mm}\ldots,h_d)\big)$.\\

 Combining the above methods, we get the $f$-vector of $P^{d,k,n}$:
\begin{equation}\label{eq4}
\textbf{f}(P^{d,k,n})=T\big(\textbf{v}+(n-k)\textbf{u},(h_{m+1},\ldots,h_d)+\textbf{0}\big)
.
\end{equation}

For an example let us return to the ordinary polytope $P^{5,7,9}$, where $\textbf{u}=(0,1,3)$ and $\textbf{v}=(1,5,10)$.

\setlength{\unitlength}{0.65mm}
\begin{picture}(100,41)(-8,1.5)

\put(75,34){$\big(1$} \put(85,34){$7$}\put(94,34){$16\big)$} \put(104.3,34){${{6}}$}

\put(70,26){$1$} \put(80,26){$8$} \put(90,26){$23$} \put(100,26){$22$} \put(110,26){${3}$} 

\put(65,18){$1$} \put(75,18){$9$} \put(85,18){$31$}
\put(95,18){$45$} \put(105,18){$25$}\put(115,18){${1}$}

\put(60,10){$\big(1$}\put(70,10){$10$} \put(79,10){$40$}
\put(89.5,10){$76$} \put(100,10){$70$}\put(110,10){$26\big)$} 

\put(30,10){$\textbf{f}(P^{d,k,n})=$}
\put(49,34){$\textbf{v}+2\textbf{u}=$}

\put(60,7.9){\line(1,0){56.2}}
\put(74,32.5){\line(1,0){36.3}}
\end{picture}

According to the formula (\ref{eq4}) and to Theorem \ref{t4} the log-concavity of $\textbf{f}(P^{d,k,n})$ depends only on the log-concavity of the vector $\textbf{v}+(n-k)\textbf{u}$, because $(h_{m+1},\ldots,h_d)$ clearly satisfies the conditions of Theorem \ref{t4}. The log-concavity of $\textbf{v}+(n-k)\textbf{u}$ is guaranteed by the following lemma, which can be verified by using the log-concavity of the rows of Pascal's triangle and using the   
recursion $\binom{n+1}{i+1}=\binom{n-1}{i-1}+2\binom{n-1}{i}+\binom{n-1}{i+1}.$

\begin{lemma}The sequence $\big\{\binom{n+1}{k}+m\binom{n-1}{k-1}\big\}_{k\in \mathbb{N}}$ is log-concave for any positive integers $m$ and $n$.
\end{lemma}

Combining the above results we get the following theorem:

\begin{theorem}The $f$-vectors of ordinary polytopes are log-concave. 
\end{theorem}


\begin{thebibliography}{9}                                                                                            
\vspace{1.2mm} 

\bibitem {bay}  \textsc{M. M. Bayer}:\ \textit{Flag vectors of multiplicial polytopes},  Elect. J. Comb., 11\vspace{0.7mm} (2004), 101-123.

\bibitem {bis}  \textsc{T. Bisztriczky}:\ \textit{Ordinary $(2m + 1)$-polytopes},  Israel J. Math., 102 (1997),\vspace{0.7mm} 101-123.

\bibitem {bjo}  \textsc{D. Bj\"orner}:\ \textit{The unimodality conjecture for convex polytopes},  Bull. Amer. Math. Soc. (N.S.) Vol. 4, Nr. 2 (1981),\vspace{0.7mm} 187-188.

\bibitem {bre} \textsc{F. Brenti}:\ \textit{Log-concave and unimodal sequences in algebra, combinatorics, and geometry: an update}, Contemporary Math., 178\vspace{0.7mm}
(1994), 71-89.

\bibitem {din} \textsc{T. N. Dinh}:\ \textit{Ordinary polytopes}, Ph.D. thesis, The University of \vspace{0.7mm}Calgary,
1999.

\bibitem {eck} \textsc{J. Eckhoff}:\ \textit{Combinatorial properties of f-vectors of convex polytopes},  Normat 54, no. 4 (2006),\vspace{0.7mm} 146-159.

\bibitem {kur}  \textsc{D. C. Kurtz}:\ \textit{A Note on Concavity Properties of Triangular Arrays of Numbers}, J. Comb. Theory, Ser. A, Vol. 13, Nr. 1 (1972),\vspace{0.7mm} 135-139.
 
\bibitem {lee} \textsc{C. W. Lee}:\ \textit{Some recent results on convex polytopes}, in: "Mathematical Developments arising from Linear Programming"  (J. C. Lagarias and M. J. Todd, eds.), Contemporary Mathematics 114 \vspace{0.7mm}(1990),
 3-19. 

\bibitem {sta} \textsc{R. Stanley}:\ \textit{The number of faces of simplicial polytopes and spheres}, Discrete Geometry and Convexity,  Ann. New York Acad. Sci. 440, edited by J.E. Goodman, et al.,\vspace{0.7mm}(1985), 212-223. 

\bibitem {wer} \textsc{A. Werner:}\ \textit{Unimodality and convexity of f-vectors of polytopes}, Preprint, TU Berlin, December 2005,\vspace{0.7mm} http://www.arXiv.org/math.CO/0512131

\bibitem {zie} \textsc{G. M. Ziegler}:\ \textit{Lectures on Polytopes}, vol.
 152 of Graduate Texts in Mathematics, Springer-Verlag, New York, \vspace{0.7mm} 1995.

\bibitem {zieg} \textsc{G. M. Ziegler}:\ \textit{Convex Polytopes: Extremal Constructions and $f$-Vector Shapes},  IAS/Park City Mathematics Series, vol. 14, (2004), 	arXiv:math/0411400v2
 
\end{thebibliography}
\end{document}